\documentclass[12pt]{article}

    \usepackage{amsmath}
    \usepackage{amsfonts}

    \newcommand{\lbl}{\label}

    \newcommand{\eq}[1]{$(\ref{#1})$}

    \newcommand{\eps}{\varepsilon}

\newcommand{\Y}{{\cal Y}}

\newcommand{\X}{{\cal X}}
\newcommand{\0}{{\bf 0}}

    \newtheorem{theo}{Theorem}
    \newtheorem{coro}{Corollary}

    \newtheorem{lemm}{Lemma}

    \def\N{\mathbb{N}}
    \def\E{\mathbb{E}}
    
    \def\0{{\bf 0}}
    
    \def\R{\mathbb{R}}

    \renewcommand{\E}{\mathbb E \,}

    \newcommand{\supp}{{\rm supp}}

    \newcommand{\eqco}{\setcounter{equation}{0}}
    \newcommand{\thco}{\setcounter{theo}{0}}
    \newcommand{\prco}{\setcounter{prop}{0}}
    \newcommand{\laco}{\setcounter{lemm}{0}}
    \newcommand{\coco}{\setcounter{coro}{0}}
    \newcommand{\cjco}{\setcounter{conj}{0}}
    
    \newcommand{\deco}{\setcounter{defn}{0}}
    \newcommand{\allco}{\eqco  \thco \prco \laco \coco \cjco \deco}
    \renewcommand{\P}{{{\cal P}}}

    \newcommand{\card}{{\rm card}}

    \newcommand{\K}{{\cal K}}

     \newcommand{\tB}{{\tilde{B}}}

    \def\bdm{\begin{displaymath}}
    \newcommand{\edm}{\end{displaymath}}
    \def\benu{\begin{enumerate}}
    \def\eenu{\end{enumerate}}
    \def\beqn{\begin{equation}}
    \def\eeqn{\end{equation}}
    \def\be{\begin{equation}}
    \def\ee{\end{equation}}
    \def\bea{\begin{eqnarray}}
    \def\eea{\end{eqnarray}}
    \newcommand{\bean}{\begin{eqnarray*}}
    \newcommand{\eean}{\end{eqnarray*}}
    \newcommand{\bear}{\begin{eqnarray}}
    \newcommand{\eear}{\end{eqnarray}}
    
    \def\R{\mathbb{R}}

    \def\qed{\hfill\hbox{${\vcenter{\vbox{
        \hrule height 0.4pt\hbox{\vrule width 0.4pt height 6pt
        \kern5pt\vrule width 0.4pt}\hrule height 0.4pt}}}$}}

\begin{document}

    \title{\bf Laws of large numbers and nearest neighbor distances}

    \author{ Mathew D. Penrose$^{*}$ and J. E. Yukich$^{**}$}

    \date{\today}
    \maketitle

\footnotetext{$~^{*}$ Research supported in part by the Alexander von
 Humboldt Foundation through a Friedrich Wilhelm Bessel Research Award}

     \footnotetext{$~^{**}$ Research supported in part by NSF grant
    DMS-0805570}


{\center{ {\Large  Dedicated to Sreenivasa Rao Jammalamadaka to mark his
65th year}} }

\begin{abstract}
We consider the sum of power weighted nearest neighbor
distances in a sample of size $n$ from a multivariate density
$f$ of possibly unbounded support. We give various criteria
guaranteeing that this sum satisfies a law of large numbers
for large $n$,
correcting some inaccuracies in the literature on the way.
Motivation comes partly from the problem of consistent
estimation of certain entropies of $f$.

\end{abstract}

\section{Introduction}

Nearest-neighbor statistics on multidimensional data are of
long-standing and continuing interest, because of their uses, for
example, in density estimation and goodness-of fit testing
\cite{BB,LQ,ZJ}, and entropy estimation \cite{BDGM, CH, KL,LPS}.
They form a multivariate analog to the one-dimensional spacings
statistics in which the work of S. R. Jammalamadaka, the dedicatee
of this paper, has featured prominently. For example, \cite{JRT}
 uses nearest neighbor balls
  to generalize the maximum spacings method to high dimensions and to establish consistency in estimation questions.

 In the present note we revisit, extend and correct
 some of the laws of large numbers
concerned with sums of power-weighted nearest-neighbor distances
that have appeared in recent papers, notably 
Penrose and Yukich \cite{PY4}, Wade \cite{Wa}, Leonenko et al. \cite{LPS}.

Fix  $d \in \N$ and $j \in \N$.
Given a finite $\X \subset \R^d$,
 and given a point $x \in \X $,
let $\card( \X)$ denote the number of elements of $\X$, and let
$D(x,\X) := D_{j}(x,\X)$ denote the Euclidean distance from $x$ to
its $j$th nearest neighbor  in the point set $\X \setminus \{x\}$,
if $\card (\X) > j$;
 set $D(x,\X):=0$ if $\card (\X) \leq j$.
Let $f$ be a probability density funticon on $\R^d$, and
let $(X_i)_{i \in \N}$ be a sequence of  independent random
$d$-vectors with common density $f$. For $n \in \N$,
let $\X_n  := \{X_1,\ldots,X_n\}$.
 Let $\alpha \in  \R$ and set
$$
S_{n,\alpha} : = \sum_{x \in \X_n } (n^{1/d}  D( x,\X_n))^\alpha
= \sum_{i=1}^n (n^{1/d}  D( X_i,\X_n))^\alpha.
$$
Certain transformations of the $S_{n,\alpha}$
have been proposed
 \cite{KL, LPS}
as  estimators for certain  entropies of the density $f$ which are
defined in terms of the integrals
$$
I_\rho(f) :=
\int_{\R^d}
f(y)^\rho dy~~~(\rho >0). 
$$

For $\rho >0$ with $\rho \neq 1$,
the
Tsallis $\rho$-entropy (or Havrda and Charv\'at $\rho$-entropy
\cite{HC}) of the density $f$ is  defined by
$H_\rho(f) := (1-  I_\rho (f))/(1 - \rho)$,
while
 the R\'enyi entropy \cite {Re} of $f$ is
defined by $H_\rho^* (f) :=  \log I_\rho(f)/(1-\rho)$.

R\'enyi and Tsallis entropies figure in various scientific
disciplines, being used in dimension estimation and the
study of nonlinear Fokker-Planck equations, fractal random walks,
parameter estimation in semi-parametric modeling, and data
compression (see \cite{CH} and \cite{LPS} for further details and
references).  

 A problem of interest is to estimate the R\'enyi and
Tsallis entropies, or equivalently, the integrals $I_\rho (f)$,
given only the sample $\{X_i\}_{i=1}^n$ and their pairwise
distances.  Let $\omega_d := \pi^{d/2} / \Gamma(1 + d/2 )$ denote the volume
of the
unit radius Euclidean ball in $d$ dimensions, and set $\gamma(d,j):= \omega_d^{-\alpha/d} \left( \frac{ \Gamma
(j + \alpha /d) }{\Gamma(j)} \right)$.  This note provides
sufficient conditions on the density $f$ establishing that $\gamma(d,j)^{-1}n^{-1}
S_{n,\alpha}$ converges to $I_{1- \alpha/d} (f)$ in $L^1$, or in
 $L^2$.
In other words, since $L^1$ convergence implies convergence of means, we
provide sufficient conditions on $f$ guaranteeing that
$\gamma(d,j)^{-1}n^{-1} S_{n,\alpha}$ is an {\em asymptotically
unbiased and consistent estimator} of  $I_{1- \alpha/d} (f)$.

\section{Results}
\lbl{secresults}

Two of our results can be stated without further ado.

\begin{theo}
\lbl{thm1} Let $\alpha > 0$.  Suppose
 the support of $f$ is a finite union of convex bounded sets
with nonempty interior, and  $f$ is
bounded away from zero and infinity on its support.
Then as $n \to \infty$ we have $L^2$ and almost sure convergence
\bea
n^{-1} S_{n,\alpha} \to \omega_d^{-\alpha/d}
\left( \frac{ \Gamma (j + \alpha /d) }{\Gamma(j)}
\right)
I_{1- \alpha/d} (f).
\lbl{conv1}
\eea
\end{theo}

\begin{theo}
\lbl{thm3} Let $q = 1$ or $q=2$. Let $\alpha \in (-d/q,0)$ and
suppose $f$ is bounded.
 Then
\eq{conv1} holds with $L^q$ convergence.
\end{theo}

For the interesting case when $\alpha > 0$ and
$f$ has unbounded support, our results require further notation.
Let $|\cdot|$ denote the Euclidean norm on $\R^d$.
For $r >0$, define the integral
$$
M_r(f) := \E[ |X_1|^{r}] = \int_{\R^d} |x|^r f(x) dx,
$$
and define the critical moment $r_c (f) \in [0, \infty]$, by
$$
r_c(f) := \sup \{ r \geq 0: M_r(f) < \infty \} .
$$
If $r < s$ and $M_{s}(f) < \infty$,
then $M_{r}(f) < \infty$. Hence
 $M_r(f) < \infty$ for $r < r_c(f)$ and
$M_{r}(f) = \infty$ for $r > r_c(f)$.

 For $k \in \N$, let  $A_k$ denote the annular
 shell centered around the origin of
$\R^d$ with inner radius $2^k$ and outer radius $2^{k+ 1}$, and let
$A_0$ be the ball centered at the origin with radius 2. For Borel
measurable $A \subset \R^d$, set $F(A) : = P[X_1 \in A] = \int_A
f(x) dx$.

We can now state the rest of our results.

\begin{theo}
\lbl{thm2} Let $q = 1$ or $q = 2$.  Let $\alpha  \in (0, d/q)$.
Suppose $I_{1-\alpha/d}(f)  < \infty$, and $r_c(f) > q \alpha
d/(d-q \alpha) $. Then  \eq{conv1} holds with $L^q$ convergence.
\end{theo}

We shall deduce from Theorem \ref{thm2},
  that when $f(x)$ decays as a power of $|x|$,
the condition $I_{1-\alpha/d}(f)<\infty$
is sufficient for $L^1$ convergence:

\begin{coro}
\lbl{thm6} 
Suppose there exists $\beta > d $ such that $f(x) =
\Theta(|x|^{-\beta})$ as $|x| \to \infty$, i.e. such that for some
finite positive $C$ we have \bea C^{-1} |x|^{-\beta} < f(x) < C
|x|^{-\beta}, ~~~~ \forall x \in \R^d , \ \ |x| \geq C. \lbl{0709c}
\eea
 Suppose also that
 $I_{1-\alpha/d}(f)<\infty$ for some $\alpha  \in (0, d)$.
 Then  \eq{conv1} holds with $L^1$ convergence.
\end{coro}

Our final result shows that in general, 
the condition $I_{1-\alpha/d}(f)<\infty$
 is {\em not} sufficient alone for $L^1$ convergence,
or even for convergence of expectations.
It can also be viewed as a partial converse to Theorem \ref{thm2}
 showing, under the additional regularity condition
\eq{reg},
that when $q = 1$ the condition $r_c(f)
> q\alpha d/(d-q\alpha)$ is close to being sharp.

\begin{theo}
\lbl{thm5} Let $ 0 < \alpha < d$. Then (i)
if $r_c(f) < \alpha d/ (d- \alpha)$, and
also for some $k_0 \in \N$ we have
\bea
0 < \inf_{k \geq k_0}  \frac{ F(A_{k}) }{ F(A_{k-1} ) }
 \leq \sup_{k \geq k_0}  \frac{ F(A_{k}) }{ F(A_{k-1} ) }
< \infty,
\lbl{reg}
\eea
then $\limsup_{n \to \infty}\E[n^{-1} S_{n, \alpha} ] = \infty$;

(ii)
for $0 < r < \alpha d/ (d- \alpha)$
    there exists a bounded continuous density function $f$ on $\R^d$
satisfying \eq{reg}, such that
 $I_{1 - \alpha/d}(f) < \infty$,
but with $r_c(f) =r$ so that
 $\limsup_{n \to \infty}\E[n^{-1} S_{n, \alpha} ] = \infty$ by part (i).
\end{theo}

\vskip.5cm The value of the limit in \eq{conv1} was
already known (see Lemma \ref{PY4lem}).
The contribution of the present paper is concerned with the
conditions under which the convergence \eq{conv1} holds; in what follows
we compare our conditions with the existing ones in the literature
and also comment on related limit results. For conditions under
which $n^{-1/2}(S_{n, \alpha} - \E S_{n, \alpha})$ is asymptotically
Gaussian, we refer to \cite{PY1, BY2, PeEJP}.

{\em Remarks.}

(i) {\em Theorem \ref{thm1}.} The condition  in Theorem \ref{thm1}
is a slight relaxation of condition C1 of
 the $L^2$ convergence  results in 
 \cite{PY4} or \cite{Wa}, which assume a polyhedral support set.
 When
the support of $f$ is the unit cube, Theorem 2.2 of \cite{JiY} gives
an alternative proof of almost sure convergence in \eq{conv1}
(we remark
that Theorem 2.2 of \cite{JiY} contains an extraneous $\E$ in the
left-hand side).
The convergence of means implied by Theorem \ref{thm1} was previously
 obtained, under some extra differentiability conditions on $f$, in 
\cite{EJS}.

(ii)  {\em Theorem \ref{thm3}. } The $L^1$ convergence of Theorem
\ref{thm3} improves upon Theorem 3.1 of \cite{LPS}, which
establishes mean convergence; the $L^2$  convergence of Theorem
\ref{thm3} is contained in Theorem 3.2 of \cite{LPS} and we include
this for completeness.

(iii) {\em Theorem \ref{thm2}.} The condition in Theorem \ref{thm2}
corrects
 the condition of the corresponding result given \cite{PY4},
where for $L^1$ convergence it is stated  that we need $r_c(f) >
d/(d-\alpha) $;
 in fact we need instead  the condition $r_c(f) > \alpha d/(d-\alpha) $.
In the proof  of Theorem \ref{thm2} below, we shall indicate the
errors in the proof in \cite{PY4} giving rise to this discrepancy.
This correction also applies to condition C2 in Theorem 2 of \cite{Wa},
the proof of which relies on the result  stated in \cite{PY4}.

(iv) {\em Theorem \ref{thm5}.} The condition \eq{reg}  holds, for
example, if $f(x)$ is a regularly varying function of $|x|$. Given
\eq{reg} and given $I_{1 - \alpha/d} < \infty$, Theorem
\ref{thm5} shows that the condition $r_c(f) \geq \alpha
d/(d-\alpha)$ is necessary for $L^1$ convergence of $n^{-1} S_{n,
\alpha}$, while
 Theorem
\ref{thm2} says that
$r_c(f) > \alpha d/ (d -\alpha)$ is
sufficient. It would be of interest to try to find
more refined necessary and sufficient conditions when
$r_c(f) = \alpha d/(d-\alpha)$.

(v) {\em General $\phi$}. For $\phi: \R^+ \to \R^+$ put $S_{n,
\phi}:= \sum_{x \in \X_n} \phi(n^{1/d}D(x, \X_n)).$  If $\phi$ has
polynomial growth of order $\alpha$, that is if there is a constant
$\alpha \in (0, \infty)$ such that $\phi(x) \leq C(1 + x^\alpha)$
for all $x \in \R^+$, then straightforward modifications of the
proofs show that under the conditions of Theorem
\ref{thm1} or Theorem \ref{thm2} we have the corresponding $L^q$ convergence
$$
n^{-1}S_{n, \phi} \to \int_{\R^d} \E [\phi(D(\0, \P_{f(x)}))] f(x)
dx,$$ where for all $\tau > 0$, $\P_{\tau}$ is a homogeneous Poisson
point process in $\R^d$ having constant intensity $\tau$, and $D(\0,
\P_{\tau})$ is the distance between the origin of $\R^d$ and its
$j$th nearest neighbor in $\P_{\tau}$.

(vi) {\em Minimal spanning trees.} Given a finite $\X \subset \R^d$
and $\phi: \R^+ \to \R^+$, let
$$
L_{\phi}(\X) := \sum_{e \in MST(\X)} \phi(|e|),
$$
where MST$(\X)$ denotes the edges in the graph of the minimal
spanning tree on $\X$.  Thus $L_{\phi}(\X)$ is the sum of the
$\phi$-weighted edge lengths in the minimal spanning tree on $\X$.
 Let $q = 1$ or $2$. If $\phi$ has polynomial growth of order
 $\alpha$, with $\alpha \in (0, d/q)$, if $I_{1 - \alpha/d}(f) <
 \infty$, and if $r_c(f) > q \alpha
d/(d-q \alpha)$ then, as may be seen by following the proof of
Theorem \ref{thm2}, the proof of Theorem 2.3(iii) of \cite{PY4} in
fact shows that as $n \to \infty$ we have
$$
L_{\phi}(\X_n) \to \frac{1}{2} \int_{\R^d} \E \left[ \sum_{e \in
MST(\0, \P_{f(x)})} \phi(|e|) \right] f(x)dx,
$$
where the convergence is in $L^q$, and where $MST(\0, \P_{f(x)})$
denotes the edges in the minimal spanning tree graph on $\0 \cup
\P_{f(x)}$ incident to $\0$, the origin of $\R^d$.  When $q = 2$,
this is new whereas for $q = 1$ and $\alpha \in (0,1)$, this
improves upon Theorem 2.3(iii) of \cite{PY4}, which requires $r_c(f)
> \max(\alpha d/(d-\alpha), d/(d-\alpha)).$

(vii) {\em Non-existence of density.}  If the $\{X_i\}_{i=1}^n$ fail
to have a density, then normalization of $S_{n, \alpha}$ may involve
exotic functions of $n$, including log periodic normalizations,  as
is the case when the $\{X_i\}_{i=1}^n$ have a Cantor distribution on
$[0,1]$; see \cite{Sh}.

(viii)  {\em Comparison with \cite{LPS}}. The convergence of
expectations corresponding to \eq{conv1}
 is given as the main conclusion in Theorem 3.1 of \cite{LPS}.  In
the case $1- \alpha/d <1$ 
of that result,
 it is claimed that this  convergence of
expectations holds without any extra conditions besides
finiteness of $I_{1 - \alpha/d}$.
 Theorem \ref{thm5} here disproves this
assertion; the argument in \cite{LPS} requires
 that convergence in distribution implies convergence of $r$th moments,
which is not in general true. 
 On the other hand, Corollary
\ref{thm6} shows that if we assume $f(x)$ decays as some power of
$|x|$ then finiteness of $I_{1 - \alpha/d}$ is indeed a sufficient condition for
convergence in $L^1$, and hence also convergence of expectations.

\section{Proofs}

\allco

This section provides the proofs of the results stated in
the preceding section.
We denote by
 $c, C,C',$ and $C''$ various strictly positive finite constants
whose values may change from line to line.
The proofs of Theorems \ref{thm1}, \ref{thm3} and \ref{thm2}
use the following result.

\begin{lemm}
\lbl{PY4lem}
Let $q \in\{1,2\}$ and $\alpha \in \R$.
Suppose for some $p > q$ that $\E[ (n^{1/d}D(X_1,\X_n))^{\alpha p}]$
is a bounded function of $n$.
Then \eq{conv1} holds with $L^q$ convergence.
\end{lemm}
{\em Proof.}
Since $D$ is a stabilizing
functional on homogeneous Poisson point processes \cite{PY4}, we
can apply Theorem 2.2 of \cite{PY4} or Theorem 2.1 of \cite{PY4}
to get $L^q$ convergence of $n^{-1}S_{n,\alpha}$
 to a limit which is  expressed as
an  integrated expectation in \cite{PY4} (see eqn (2.15) of
\cite{PY4}). It was shown in
\cite{Wa} that this limit is equal to the right hand side of
\eq{conv1} (and this is also consistent with the limiting constant
in \cite{EJS}). $\qed$ \\

 {\em Proof of Theorem \ref{thm1}.} Recall that we assume the support of $f$, namely $\supp (f) := \{x
\in \R^d: f(x) >0\}$,  is a finite union of bounded  convex sets
with nonempty interior, here denoted $B_1,\ldots,B_m$.
 Set $\lambda : = \sup
\{|x -y|: x \in \supp (f) , y \in \supp (f)\}$, the diameter of
the support of $f$. By assumption,
 $\lambda < \infty$. Also we assert that
there is a constant $c > 0 $ such that for $r \in (0,
\lambda]$,
 \bea
 F(B_r(x)  )  \geq  c  r^d, ~~~~   \forall x \in  \supp(f),
\lbl{0629a} \eea  
To see this, take $\delta_1 >0$ such that for $1 \leq i \leq m$
there is a ball $B_i^-$ of radius $\delta_1$ contained in $B_i$.
There is a constant $\delta_2 >0$ such that for $ 1 \leq i \leq m$,
if $x \in B_i$, and
 $r \leq \delta_1$, then
the  intersection of the ball of radius $r$ centered at $x$ with the
 convex hull of the union of $B_i^-$ and $x$
has volume at least $\delta_2 r^d$. This region is contained in
$B_i$ and \eq{0629a} follows for $r \in (0, \delta_1]$. But then
(with a different choice of $c$) \eq{0629a} follows for $r \leq
\lambda$. Hence, for $0< t \leq  \lambda n^{1/d} $  and
with $B(x;r)$ denoting the Euclidean ball of radius $r$
centered at $x$,
 \bean P [ n^{1/d}
D(X_1,\X_n) > t] \leq \sup_{x \in \supp(f)} P[ \card (\X_{n-1} \cap
B(x;n^{-1/d} t ) ) < j]
\\
\leq \sum_{i=0}^{j-1} { n-1 \choose i } ( c n^{-1} t^d)^i ( 1-
c n^{-1}t^d)^{n-1-i}
\\
\leq C \sum_{i=0}^{j-1}
t^{id} \exp ( - c n^{-1} t^d (n-1-i) )
\\
\leq C' 
( 1 + t^{(j-1)d}) \exp   ( - c t^d)
\leq C'' \exp ( - (c/2) t^d ). \eean 
Moreover this probability
is clearly zero for $t > \lambda n^{1/d} $. Hence, 
for $\alpha > 0$ and $p > 2$,
\bean 
\E[
(n^{1/d} D( X_1,\X_n))^{\alpha p} ] = \int_0^\infty P[ n^{1/d} D(
X_1,\X_n) > u^{1/(\alpha p)} ] du
\\
 \leq  C \int_0^\infty \exp( - (c/2) u^{d/(\alpha p) } ) du
\eean 
which is finite and does not depend on $n$. Therefore we can
apply 
Lemma \ref{PY4lem} to get the $L^2$ convergence
\eq{conv1}.

For almost sure convergence, we apply Theorem 2.2 of \cite{PeBer},
where here the test function considered in that result (and denoted
$f$ there, not to be confused with the notation $f$ as used here)
 is the identity function.  It is
well known (see \cite{BB}, or Lemma 8.4 of \cite{Yubk})  that there
is a constant $C:=C(d)$ such that for any finite $\X \subset \R^d$,
any point $x \in \X$  is the $j$th nearest neighbor of at most $C$
other points of $\X$. Therefore adding one point to a set $\X$
within the bounded region $\supp(f)$ changes the sum of the
power-weighted $j$th nearest neighbor distances by at most a constant.
Therefore (2.9) of \cite{PeBer} holds here (with $\beta =1$ and
$p'=5$ say), and the almost sure convergence follows by Theorem 2.2
of \cite{PeBer}. $\qed$ \\

{\em Proof of Theorem \ref{thm3}.} The proof depends on the
following lemma. Recall that $(X_i)_{i \geq 1}$ are i.i.d. with density $f$.
Given $X_1$, let $V_n$ denote the volume of the $d$-dimensional
ball centered at $n^{1/d}X_1$ whose radius equals the distance to
the $j$th nearest point in $n^{1/d}(\X_n \setminus X_1)$,
where for $r > 0$ and $\X \subset \R^d$ we write
$r \X$ for $\{rx: x \in \X\}$. For all
$x \in \R^d$, for all $n = 2,3,...$ and for all $v \in (0,
\infty)$ let \be \label{cdf}
 F_{n,x}(v):= P[V_n \leq v | X_1 = x].
\ee

\begin{lemm} \label{intparts1} If $f$ is bounded
 and $\delta \in (0,1)$, then
$$
\sup_n \E V_n^{-\delta}= \sup_n \int_{\R^d} \int_0^{\infty}
v^{-\delta} dF_{n,x}(v) f(x) dx < \infty.
$$
\end{lemm}

\vskip.5cm

{\em Proof of Lemma \ref{intparts1}.}  Since $\int_0^{\infty}
v^{-p} dF(v) = p \int_0^{\infty} v^{-p-1} F(v)dv$ for any $p\in
(0,1)$ whenever both integrals exist (see e.g. Lemma 1 on p. 150
of \cite{Fe}), 
we have for all $x \in \R^d$
\bean
 \int_0^{\infty} v^{-\delta} dF_{n,x}(v) = \delta  \int_0^{\infty} v^{-\delta-1} F_{n,x}(v)dv
\\
\leq \int_0^{1} v^{-\delta-1} F_{n,x}(v)dv + \delta \int_1^{\infty}
v^{-\delta-1} dv = \int_0^{1} v^{-\delta-1} F_{n,x}(v)dv + 1.
\eean
With $\tB_v(x)$ denoting the ball of volume $v$ around $x$,
 for all $v \in (0,1)$ we have
\bea
F_{n,x}(v) = P[V_n \leq v | X_1 = x] = 1 -
P[\text{card}(n^{1/d}\X_{n-1} \cap \tB_v(n^{1/d}x)) < j]
\nonumber \\
\leq  1 - P[\text{card}(n^{1/d}\X_{n-1} \cap \tB_v(n^{1/d}x)) = 0]
\nonumber  \\
 \label{JY1} = 1 - \left(1 - \int_{\tB_{v/n}(x)}
f(z)dz\right)^{n-1}. 
\eea
 Since $f$ is assumed bounded we have
$$
F_{n,x}(v) \leq 1 - \exp\left( (n - 1) \log (1 - \|f\|_{\infty}
v/n) \right).
$$
When $n$ is large enough, 
for all $v \in (0,1)$ 
we have 
 $(n - 1)\log (1 - \|f\|_{\infty} v/n) \geq -2 \|f\|_{\infty} v$, and so
for all $x \in \R^d$
$$
F_{n,x}(v) \leq 1 - \exp(-2\|f\|_{\infty}v) \leq 2 \|f\|_\infty v.
$$
Hence
for all $n$ large enough and all $x$ we have  $\int_0^{1}
v^{-\delta-1} F_{n,x}(v)dv \leq 2 \|f\|_\infty \int_0^{1} v^{-\delta-1} v dv 
,$ demonstrating Lemma
\ref{intparts1}.  \qed

\vskip.5cm

Now to prove Theorem \ref{thm3}, we choose $p > q$ such that $-1 <
\alpha p/d < 0$ and invoke Lemma \ref{intparts1} to conclude
$\sup_n \E[V_n^{\alpha p/d}] < \infty$.  We now apply Lemma
\ref{PY4lem}
 to complete the proof of $L^q$ convergence.  \qed \\

The proof of Theorem \ref{thm2} uses the following  lemma.
Recall from Section \ref{secresults} the
 definition of the regions $A_k, k \geq 0$.

\begin{lemm}
\lbl{Claim} Let $0 < s < d$. If
 $r_c(f) > sd/(d-s)$, then $\sum_{k=1}^\infty 2^{ks}
(F(A_k)) 
^{(d-s)/d} < \infty$.
\end{lemm}

{\em Proof.} We modify some of the arguments on page 85 of \cite{Yubk}.
For all $\eps > 0$,
by H\"older's inequality
we have
\bean
 \sum_k 2^{ks} (F[A_k])^{(d-s)/d} = \sum_k 2^{-\eps k s} (F[A_k])^{(d-s)/d} 2^{( 1 +
\eps)ks}
\\
\leq \left( \sum_k (2^{-\eps k s})^{d/s} \right)^{s/d} \left(
  \sum_k F[A_k] (2^{( 1 + \eps)ks})^{d/(d-s)}
  \right)^{(d-s)/d}
\\
 \leq C(\eps, s) \left( \sum_k \int_{A_k} |x|^{(1 +
  \eps)sd/(d-s)} f(x) dx \right)^{(d-s)/d}
\eean
  which, for $\eps$ small enough,  is finite by hypothesis.  \qed \\

{\em Proof of Theorem \ref{thm2}.} We follow the proof in
\cite{PY4}, but correct it in some places and give more details in
others.
We aim to use Lemma \ref{PY4lem}.
Since we assume $0< \alpha <d/q$, we can take $p > q$ with
$\alpha p < d$.
  Clearly
 \bea
\E[ (n^{1/d}D(X_1,\X_n))^{\alpha p}] =
n^{\alpha p/d -1 }
\E \left[ \sum_{i=1}^n D(X_i,\X_n)^{\alpha p}  \right]
=
n^{\alpha p/d -1 }
\E [ L^{\alpha p} (\X_n)],
\lbl{090622b}
\eea
where
for any finite point set $\X \subset \R^d$,
and any  $b >0$, we write
$L^b(\X) $ for $\sum_{x \in \X} D(x,\X)^b$
(and set $L^b(\emptyset) := 0$).
Note that for some finite constant $C=C(d,j)$ the functional
$\X \mapsto L^b(\X)$ satisfies
the simple subadditivity relation
\bea
L^b(\X \cup \Y ) \leq L^b(\X) + L^b (\Y) + C t^b
\lbl{subadd}
\eea
for all $t >0 $ and all finite $\X$ and $\Y$
contained in $[0,t]^d$ (cf. (2.2) of \cite{Yubk}).

As in (7.21) of \cite{Yubk} or (2.21) of \cite{PY4}
we have that
\bea
L^{\alpha p}(\X_n) \leq \left( 
\sum_{k =0} ^\infty
L^{\alpha p}( \X_n \cap A_k ) \right) + C(p) \max_{ 1 \leq i \leq n} |X_i|^{\alpha p}.
\lbl{090622a}
\eea
In the last sentence of the proof of Theorem 2.4 of \cite{PY4} it
is asserted that
 the last term in \eq{090622a}  is not needed, based on
a further assertion that one can take $C=0$ in \eq{subadd}  here,
but these assertions are incorrect. For example,
if $\card(\Y) \leq j$ then $L^b(\Y) =0$ but $L^b(\X \cup \Y)$
could  be strictly greater than $L^b(\X)$.
Similarly,  if $\card (\X_n \cap A_k ) \leq j $
 then the term in \eq{090622a} from that $k$ is
 zero but the corresponding contribution to
 the left side  of \eq{090622a} is non-zero.

Combining \eq{090622a} with \eq{090622b} yields 
\bea
\E[
(n^{1/d}D(X_1,\X_n))^{\alpha p}] \leq n^{(\alpha p -d)/d} \E \left[
\sum_{k} L^{\alpha p}(\X_n \cap A_k) \right]
\nonumber \\
  \ \ \ \ \ \ \ \ \ \  \ \
+ C(p) \E [ n^{(\alpha p -d)/d} \max_i  |X_i| ^{\alpha p} ].
\lbl{090622c} 
\eea
 By Jensen's inequality and the growth bounds
 $L^{\alpha p}(\X) \leq C
(\text{diam} \X)^{\alpha p} (\text{card} (\X))^{(d-\alpha p)/d}$
(see Lemma 3.3 of \cite{Yubk}),
we can  bound the first term in the right hand side of \eq{090622c}
by 
\bea C \sum_k  2^{k \alpha p} (F[A_k])^{(d - \alpha p )/d}.
\lbl{090622d} \eea
Recall that we are assuming $0< \alpha < d/q$ and also
 $r_c(f) > qd \alpha/(d-q \alpha)$
(the last assumption  did not feature
 in \cite{PY4}, but in fact we do need it).
Let $p > q$ be chosen so that
 $r_c(f) > d \alpha p/(d- \alpha p)$ as well as
 $\alpha p < d$.
Setting $s =\alpha p$ in Lemma \ref{Claim}, we get that
 the expression \eq{090622d} is finite.
Thus the first term in the right hand side of \eq{090622c}
is bounded by a constant independent of $n$.

The second term in the right hand side of \eq{090622c}
is bounded by
\bean
 & C(p) & \left(
 \int_0^{1}
P \left[ \max_{1 \leq i \leq n} |X_i|^{\alpha p}
\geq t n^{(d-\alpha p)/d} \right] dt
+ \int_1^\infty
P \left[ \max_{1 \leq i \leq n} |X_i|^{\alpha p}
\geq t n^{(d-\alpha p)/d} \right] dt \right)
\\
\leq  &C(p) & \left( 1 + n  \int_1^\infty P[ |X_1|^{\alpha p d /(d - \alpha p)} \geq t^{d/(d- \alpha p)} n ] dt
\right) .
\eean
By  Markov's inequality together with the assumption
$r_c(f) >  d \alpha p /(d- \alpha p)$, this last integral is
bounded by a constant independent of $n$.

Therefore the expression \eq{090622c} is bounded independently of
$n$,  so we can apply Lemma \ref{PY4lem}
to
get the $L^q$ convergence in \eq{conv1}.
$\qed$ \\


{\em Proof of Corollary \ref{thm6}.} Suppose for some
 $\beta > d$  that
$f(x) = \Theta(|x|^{-\beta})$ as $|x| \to \infty$. Then it is easily
verified that given $\alpha \in (0,d)$, the condition
$I_{1-\alpha/d}(f) < \infty$ implies that $-\beta(1-\alpha/d) +d <0$
and hence $\beta > d^2(d-\alpha)^{-1}$. Moreover, it is also easily
checked that $r_c (f) = \beta -d$ so that if $\beta >
d^2(d-\alpha)^{-1}$ then $r_c (f) > d \alpha /(d - \alpha)$.

Therefore, if $I_{1-\alpha/d}(f) < \infty$ we can apply the case
$q=1$ of Theorem \ref{thm2} to get \eq{conv1} with $L^1$
convergence.
$\qed$ \\

The proof of Theorem \ref{thm2} shows that  \be \label{manifolds}
\E[ (n^{1/d} D(X_1, \X_n))^{\eps}] < \infty, \ee if $\eps >
0$ is such that $\eps d/(d - \eps) < r_c(f)$.  The proof of
Theorem \ref{thm5}, given below, shows that the condition $\eps
d/(d - \eps) < r_c(f)$ cannot be dropped in general. \vskip.5cm

{\em Proof of Theorem \ref{thm5}.} Let $0 < \alpha < d$. Suppose
that  $ r_c(f) <  \alpha d / (d- \alpha)$, and \eq{reg} holds for
some $k_0 \in \N$.
Choose $r,s $ such that
  $r_c(f)  < r < s <  \alpha d / (d- \alpha)$.
Then $M_r(f) = \infty$, so $\sum_k 2^{rk} F(A_k) = \infty$ and
therefore there is an infinite subsequence $\K$ of $\N$ such that
\bea 2^{sk } F(A_k)  \geq 1,  ~~~ k \in \K. \lbl{calK} \eea Indeed,
if no such $\K$ existed, then for all but finitely many $k$ we would
have $2^{rk} F(A_k) \leq 2^{(r-s)k}$ which is summable in $k$.

Given $k \in \N$, and set $ n(k) = \lceil (F(A_k))^{-1} \rceil$, the
smallest integer not less than
 $(F(A_k))^{-1}$. Let $E_k$ be the event that $X_1 \in A_k$ but
$X_i \notin A_{k-1} \cup A_k \cup A_{k+1}$ for $2 \leq i \leq n(k)$.
Then by the condition \eq{reg},
  there is a strictly positive constant $c$, independent of $k$,
such that for $k \geq k_0$ we have
$$
P[E_k]   = F(A_k) (1-F(A_{k-1} \cup A_k \cup A_{k+1}))^{n(k)-1}
\geq c F(A_k) .
$$
If $E_k $ occurs then $D(X_1,\X_{n(k)}) \geq 2^{k-1}$, so for $n =
n(k)$ we have
 (for a different constant $c$) that
\bean \E [ n^{-1} S_{n,\alpha} ] = \E[ (n^{1/d} D(X_1,\X_n))^\alpha]
\geq
 n^{\alpha/d} \E[ D(X_1,\X_n)^\alpha {\bf 1}_{E_k}]
\\
\geq c n^{\alpha/d}  F(A_k) 2^{k \alpha} \geq  c (F(A_k))^{1-
\alpha /d} 2^{k \alpha}.
\eean By \eq{calK}, for $k \in \K$ this lower bound is at least a
constant times $   2^{k(\alpha - s (d-\alpha)/d)}$, and
therefore tends to infinity as $k \to \infty$ through the sequence
$\K$, concluding the proof of part (i).

 For part (ii), for each $k \geq 2$ choose, in an arbitrary way,
 a unit radius ball $B_k$ that is 
 contained in $A_k$.
 Given $r \in (0, \alpha d/(d- \alpha) )$,
 consider the density function $f$ with $f(x) = C 2^{-rk}$
 for $x \in B_k, k \geq 2$, and with $f(x) =0$
 for $x \in \R^d \setminus \cup_{k=2}^\infty B_k$; here the normalizing
 constant $C$ is chosen  to make $f$ a probability
density function.
This gives $F(A_k) = C \omega_d 2^{-rk}$
for  each $k\geq 2$; it is easy to see that this $f$ has $r_c(f)
=r$, and that \eq{reg} holds with $k_0=3$. Also,
for any $\rho >0$ we have
$
I_{\rho} (f) =  \omega_d C^\rho \sum_{k \geq 2} 2^{-r \rho k}
$
which is finite, so in particular $I_{1-\alpha/d}  < \infty$.
This choice of $f$ is bounded but not continuous, but can easily be modified
to a continuous density with the same properties, for example
by modifying $f$ in an annulus near the boundary of each ball
$B_k$ so as to make it continuous, and then adjusting the normalizing
constant $C$ accordingly.
$\qed$

Mathew D. Penrose, Department of Mathematical Sciences, University
of Bath, Bath BA2 7AY, United Kingdom:
{\texttt m.d.penrose@bath.ac.uk } \vskip.5cm

J. E. Yukich, Department of Mathematics, Lehigh University,
Bethlehem PA 18015:
\\
{\texttt joseph.yukich@lehigh.edu}

\end{document}